\documentclass{article}[11]

\usepackage{graphics,latexsym}
\usepackage{graphicx}
\usepackage{amsmath, amssymb,amsthm}

\usepackage{color}

\newcommand{\R}{\mathbb{R}}
\begin{document}

\title{Isospectral domains for discrete elliptic operators\small\footnote{
This work has been partially  supported by GNCS-INDAM (Gruppo Nazionale
Calcolo Scientifico - Istituto Nazionale di Alta Matematica).
}}

\author{Lorella Fatone$^{1}$, Daniele Funaro$^{2}$}
\date{}
\maketitle

\centerline{$^{1}$\small Dipartimento di Matematica  e Informatica }
\centerline{\small Universit\`a di Camerino, Via Madonna delle Carceri 9, 62032
Camerino (Italy)} 
\medskip

\centerline{$^{2}$\small Dipartimento di Fisica, Informatica e Matematica }
\centerline{\small Universit\`a di Modena e Reggio Emilia, Via Campi 213/B, 41125
Modena (Italy)} 
\medskip

\medskip

\begin{abstract}
Concerning the Laplace operator with homogeneous Dirichlet
boundary conditions, the classical notion of {\sl isospectrality} assumes
that two domains are related when they give rise to the same spectrum. In two dimensions, non isometric,
isospectral domains exist. It is not known however if all the eigenvalues relative
to a specific domain can be preserved under suitable continuous deformation of its geometry.
We show that this is possible when the 2D Laplacian is replaced by a finite dimensional
version and the geometry is modified by respecting certain constraints. The
analysis is carried out in a very small finite dimensional space, but it can
be extended to more accurate finite-dimensional representations of the 2D Laplacian, with an increase of
computational complexity. The aim of this paper is to introduce the preliminary
steps in view of more serious generalizations. 
\end{abstract}

%
%

\section{Introduction}\label{sec1}

Consider the Laplace problem in 2D with homogeneous
Dirichlet boundary conditions defined in an open set with regular boundary
(see \cite{grebe} for a general overview).
It is known that there are distinct domains (non isometric) such that
all the infinite eigenvalues of the Laplace operator coincide (see, e.g., Fig.~\ref{fig1}). 
For this reason,
these are called isospectral domains. It can be shown that two isospectral
domains have the same area. It is not known however
if it is possible to connect with continuity two isospectral domains
through a sequence of domains, by preserving the whole spectrum.

Partial answers can be given by working in a finite dimensional environment.
Here, we take a suitable approximation of the Laplace operator corresponding
to a negative-definite matrix. By varying the domain,
we are interested to detect those deformations that preserve the entire set
of eigenvalues (that are now in finite number). At the same time, not
all the possible deformations are allowed, but only those belonging
to a finite dimensional space of parameters. The problem turns out to
be far from easy. Indeed, already at dimension 4, things get rather involved.
The question examined here concerns with the deformation
of quadrilateral domains, with the aim of preserving the eigenvalues
of discrete operators obtained from collocation of the Laplace problem
using polynomials of degree 3 in each variable. By imposing Dirichlet  homogeneous boundary
conditions, the discretization matrix ends up to be only of dimension $4\times 4$. 
The vertices of the domains are then suitably moved by maintaining the 
magnitude of the four corresponding eigenvalues. The results show that,
at least in these simplified circumstances, families of isospectral domains
exist and can be connected by continuous transformations. 

The most straightforward
(but extremely expensive) approach is to try all the possible allowed
configurations and sort them by comparing the so obtained spectrum. 
Upgraded versions consist in moving the vertices along curves, whose tangent
is obtained as an application of the {\sl Implicit Function} theorem due to U. Dini (see, for example, \cite{rudin}).
 Here, one computes
with the help of symbolic manipulation the partial derivatives, with respect to the 
various parameters,  of the coefficients of the 
characteristic polynomial of the matrix representing the discrete operator. 
Different, more or less efficient, approaches have been tested.
By the way, despite its simple formulation, the problem looks rather complex and the
extension to higher dimensions or to more complicated domains looks
at the moment quite unrealistic.

\section{Preliminary settings}\label{sec2}

For the convenience of the reader we  briefly review some results about the   classical eigenvalue problem for 
the Laplace operator:
 $\Delta= \frac{\partial^{2}}{\partial x^{2}}+\frac{\partial^{2}}{\partial y^{2}}$, 
in an open set $\Omega$,  when  homogeneous Dirichlet boundary conditions are  imposed on the 
boundary $\partial \Omega$. The  problem is formulated as follows:
\begin{eqnarray}
& -\Delta u(x,y)=\lambda u(x,y), \qquad  \ &( x, y)  \ \in  \ \Omega, \label{lapl1} \\
& u(x,y)=0 ,\qquad  \  &( x, y)  \ \in   \ \partial \Omega. \label{laplBC1}
\end{eqnarray}
  It is known (see, e.g.: \cite{grebe}) that the spectrum of minus the Laplace operator 
 is discrete, that the  eigenvalues are non negative and can be ordered in  
ascending order to form a divergent sequence:
\begin{eqnarray} \label{aut1}
&& 0 < \lambda_{1} < \lambda_{2} \le \lambda_{3} \le \ldots  \nearrow  \infty ,
\end{eqnarray}
with possible multiplicities. When  $\Omega=Q=]0,1[\times ]0,1[$ the eigenvalues are:
\begin{eqnarray} \label{aut1Q}
&&  \lambda  =\pi^{2} (m^2+n^2), \quad m,n=1,2,3,4, \ldots.
\end{eqnarray}

The Weyl's law establishes an estimate of the
$m$-th eigenvalue in terms of the area $\mu_2 (\Omega )$ of the domain $\Omega$
where the eigenvalue problem is defined: 
\begin{eqnarray}
\frac{\lambda_m}{m} \rightarrow \frac{4\pi }{\mu_2 (\Omega )},   \quad {\rm for } \ m \rightarrow \infty .
 \label{weyl}
\end{eqnarray}
This relation leads us to the following consequence. For a given $\lambda$, one
denotes by $N(\lambda )$ the number of eigenvalues smaller than $\lambda$.
Then, a Weyl's theorem states that:
\begin{eqnarray}
N( \lambda) \approx \frac{\mu_2 (\Omega )}{4\pi}\lambda.  
 \label{weyl2}
\end{eqnarray}
It turns out that, if two distinct domains produce the same set of
eigenvalues (i.e., they are isospectral), they must have the same area.
A refinement of property  (\ref{weyl2}) brings to the relation:
\begin{eqnarray}
N( \lambda) \approx \frac{\mu_2 (\Omega )}{4\pi}\lambda -\frac{\mu_1 (\partial\Omega )}
{4\pi}\sqrt{\lambda} ,
 \label{weyl3}
\end{eqnarray}
known as Weyl's conjecture, which also involves the length of the perimeter $\partial\Omega$, i.e.  $\mu_1 (\partial\Omega )$.
This means that, provided the conjecture is verified, isospectral domains
also share the same perimeter.

\begin{figure}
\centerline{\includegraphics[height=5cm]{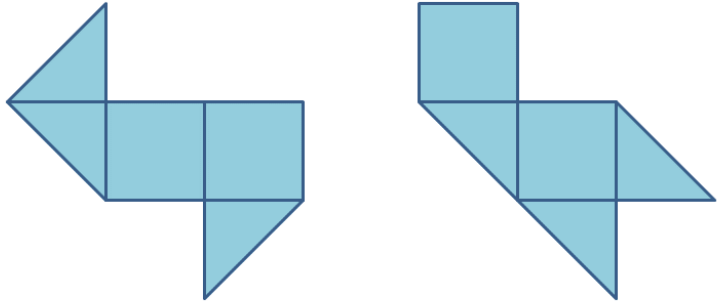} }
  \caption{Two  isospectral domains in $\R^{2}$. They have the same area 
and the same perimeter. Moreover, they have the same {\sl sound}, i.e.
they share an identical sequence of Laplace eigenvalues. }
 \label{fig1}
  \end{figure}

The existence of (non isometric) isospectral domains has been established quite
recently. The problem was firstly formulated in \cite{kac}, where the author
was wondering if the eigenspectrum of the Laplacian was sufficient to detect
the shape of the domain (put in other words: if one can {\sl hear} the shape of a drum). 
In the 2D plane, the negative answer appeared in \cite{gordon},
where distinct domains (see, for example,  Fig.~\ref{fig1}) exhibiting  
identical sets of eigenvalues, were proposed. Preliminary results in this direction were 
investigated in \cite{sunada}.
Successive examples have been discussed for instance in \cite{buser} and
\cite{chapman}. Nowadays, wide classes of isospectral domains are available.
For the sake of brevity, we address the reader to the specialized literature for more insight.

Some facts are still not known however, as for instance the existence of 2D isospectral
domains of convex type, or the possibility to vary with continuity the shape
of a given domain maintaining at the same time its entire spectrum. This last 
property is the one we are going to investigate in this paper, when the Laplace operator
is substituted by a very rough finite-dimensional version. In fact, in this paper we want to see 
if, under suitable simplified  hypotheses, it is possible to connect with continuity two isospectral domains
through a sequence of domains, by preserving the whole spectrum.

Our first goal is to define a general set of quadrilateral domains. Afterwords, 
starting from a given member of the family, we will try to find 
other members that are isospectral to it. The notion of isospectrality is
decided according to a finite-dimensional elliptic operator that is
going to be introduced in Sect.~\ref{sec3}. First of all, we need to exclude as much
as possible, the chance of having isometric domains in the family, i.e,
quadrilaterals that can be related through elementary operations, such 
as translation, rotation and symmetry. Of course, two isometric domains
are automatically isospectral and we would like to avoid such trivial connections.

 \begin{figure}
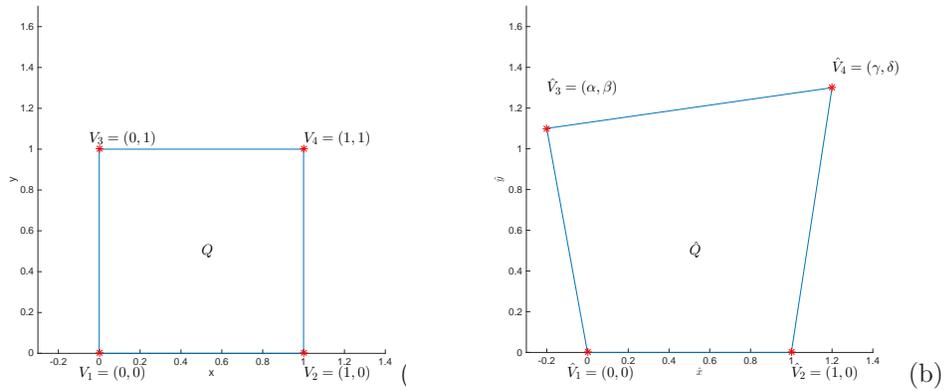

\centerline{\includegraphics[height=5cm]{fig2AQuadrato.eps} (a)\hskip.8truecm\includegraphics[height=5cm]{fig2BQSbilenco.eps} (b)}
  \caption{(a) The unit square $Q$. (b) The generic quadrilateral $ \hat{Q}$.}
 \label{fig2}
  \end{figure}

From now on we assume that one of the sides of our quadrilaterals is ``nailed" to the $x$-axis of the plane.  
Moreover in the future, we 
shall avoid those cases in which a couple of sides may intersect at some internal point.

Let $n$ be a positive integer, $ \R$ be  the set of  real numbers, $ \R^{n}$ be  the $n$-dimensional real Euclidean space.
Let $Q \subset \R^{2}$ be the unit square in the space $\R^{2}$, i.e., 
$Q$ is the quadrilateral  of vertices  $V_{1}=(0,0)$, $V_{2}=(1,0)$, $V_{3}=(0,1)$, $V_{4}=(1,1)$ 
(see Fig.~\ref{fig2}). 
Given the real parameters $ \alpha ,\beta, \gamma , \delta$, let 
 $\hat{Q} \subset \R^{2}$  be the generic quadrilateral of our family, having
vertices  $\hat{V}_{1}=(0,0)$, $\hat{V}_{2}=(1,0)$, $\hat{V}_{3}=(\alpha ,\beta)$, $\hat{V}_{4}=(\gamma ,\delta)$
(see again Fig.~\ref{fig2}).  The domains $Q$ and $ \hat{Q}$  are open.

In the plane with coordinates $(\hat{x},\hat{y})$,  we now focus our attention on the  eigenvalue problem 
for the Laplace operator defined in $\Omega =\hat{Q}$,  with homogeneous Dirichlet boundary conditions on the piecewise 
smooth  boundary $\partial \hat{Q}$. Translating into formulas, we have:
\begin{eqnarray}
& -{\Delta} \hat{u} (\hat x, \hat y) ={\lambda} \hat{u} (\hat x, \hat y),  \qquad  \ & (\hat x, \hat y) \ \in \ \hat{Q}, \label{lapl2} \\
&\hat{u} (\hat x, \hat y)=0 ,  \qquad  \ & (\hat x, \hat y) \ \in    \ \partial \hat{Q}, \label{laplBC2}
\end{eqnarray}
where $\displaystyle{\Delta}= \frac{\partial^{2}}{\partial \hat{x}^{2}}+\frac{\partial^{2}}
{\partial \hat{y}^{2}}$. 

For convenience, let us  map  problem  (\ref{lapl2}), (\ref{laplBC2})   into the following
modified version, defined in the square $Q$: 
\begin{eqnarray}
& L  u ( x, y)  ={\lambda}  u ( x, y) , \qquad  \ &( x, y)  \ \in  \   Q, \label{lapl3} \\
&  u ( x, y) =0 , \qquad  \ &( x, y)  \ \in  \  \partial  Q, \label{laplBC3}
\end{eqnarray}
where $L$ turns out to be a suitable positive-definite elliptic operator that we are going to define
here below.
To this scope let us examine   the transformation $\hat x=\theta_1 (x,y)$, $\hat y=\theta_2 (x,y)$ that 
allows us to bring the operator $\Delta$ into $L$. 
First of all let us  transform a general quadrilateral $ \hat{Q}$
into the reference square $Q$. We use a classical invertible mapping $\theta: Q \rightarrow \hat{Q}$ 
consisting of polynomials of degree one in each variable. This relates the ordered points
$V_{1}$,  $V_{2}$, $V_{3}$, $V_{4}$ of  $Q$ with the ordered points
$\hat{V}_{1}$,  $\hat{V}_{2}$, $\hat{V}_{3}$, $\hat{V}_{4}$ of $ \hat{Q}$ (see Fig.~\ref{fig2}). 
The two components of the transformation  $\theta=(\theta_{1},\theta_{2})$ are given by:
\begin{eqnarray}
&& \theta_{1}:   \hskip1.2truecm \hat{x}=x+\alpha y +(\gamma-1-\alpha) x y, \label{trasf1} \\
&& \theta_{2}:   \hskip1.2truecm \hat{y}=\beta y +(\delta-\beta) x y.\label{trasf2}
\end{eqnarray}
Thus, a function $\hat u$  defined in $\hat{Q}$ is associated with the function $u=\hat u(\theta)$ defined in $Q$.
By applying the change of variables to the Laplace operator $\Delta$ we arrive at
the eigenvalue problem (\ref{lapl3}), (\ref{laplBC3})  where the operator $L$ is  defined as follows (see  \cite{miolibro2}):
\begin{eqnarray}\label{defL} 
&&\hskip-1truecm L =f_{1} \frac{\partial^{2}}{\partial x^{2}}  + f_{2} \frac{\partial^{2}}{\partial x \partial y}  +  f_{3} \frac{\partial^{2}}{\partial y^{2}}  +
f_{4} \frac{\partial }{\partial x}  + f_{5} \frac{\partial }{\partial y}, 
\end{eqnarray}
and  the coefficients of $L$ in (\ref{defL}) are given by:
\begin{eqnarray}\label{fi}
&& \hskip-0.7truecm f_{1}= -\frac{1}{\sigma^{2}} \left[   \left(  \frac{\partial\theta_{1}}{\partial y}     \right)^{2}+
		   \left(  \frac{\partial \theta_{2}}{\partial y}     \right)^{2}
                  \right], \quad 
 f_{2}= \frac{2}{\sigma^{2}} \left[     \frac{\partial\theta_{1}}{\partial x}   \frac{\partial\theta_{1}}{\partial y}   +
		     \frac{\partial\theta_{2}}{\partial x}   \frac{\partial\theta_{2}}{\partial y} 
                  \right], \nonumber\\[3mm]
&&\hskip-0.7truecm f_{3}= -\frac{1}{\sigma^{2}} \left[   \left(  \frac{\partial\theta_{1}}{\partial x}     \right)^{2}+
		   \left(  \frac{\partial \theta_{2}}{\partial x}     \right)^{2}
                  \right],\quad
 f_{4}= \frac{ f_{2}}{\sigma} \left[     \frac{\partial\theta_{1}}{\partial y}   \frac{\partial^{2}\theta_{2}}{\partial x \partial y}   -
		     \frac{\partial\theta_{2}}{\partial y}   \frac{\partial^{2}\theta_{1}}{\partial x \partial y} 
                  \right],\\[3mm]
&& \hskip-0.7truecm f_{5}= \frac{ f_{2}}{\sigma} \left[     \frac{\partial\theta_{2}}{\partial x}   \frac{\partial^{2}\theta_{1}}{\partial x \partial y}   -
		     \frac{\partial\theta_{1}}{\partial x}   \frac{\partial^{2}\theta_{2}}{\partial x \partial y} 
                  \right], \nonumber
\end{eqnarray}
where
\begin{eqnarray}
&& \sigma=  \frac{\partial\theta_{1}}{\partial x}     \frac{\partial\theta_{2}}{\partial y}  -  \frac{\partial\theta_{1}}{\partial y}    \frac{\partial\theta_{2}}{\partial x},  
\end{eqnarray}
is the determinant of the Jacobian of $\theta$.

In the future, instead of solving the  eigenvalue problem directly on $\hat{Q}$, we will find more convenient to approach
the transformed  eigenvalue problem defined in $Q$.  That is, instead of approaching problem  (\ref{lapl2}), (\ref{laplBC2}) we consider problem
  (\ref{lapl3}), (\ref{laplBC3}) where $L$ is  defined in   (\ref{defL}). It is straightforward to
check that $L=\Delta$ when $\hat{Q}=Q$.

As explained with more details  in Sect.~\ref{sec4}, it is also convenient
to introduce a new parameter $c>0$. Through the homothety centered at the origin  $(0,0)$, we associate
a generic point $\hat p$ of the plane $(\hat x, \hat y)$ to the point $\sqrt{c}  \ \hat p$. 
In this way, any given eigenvalue $\lambda$ of $\Delta$ in the domain $\hat{Q}$ takes
the form $\lambda /c$ in the new domain $\sqrt{c} \ \hat{Q}$. By this trick, without too much effort, we can
enlarge the set of quadrilaterals at our disposition.

\section{The discrete operator}\label{sec3}

In the referring domain $Q$ we now build a discrete elliptic operator
corresponding to a very rough approximation of the operator $L$ defined
in (\ref{defL}). Indeed,
we do not want to handle too many eigenvalues. The reason is that the
family of quadrilateral domains introduced in the previous section only
depends on five degrees of freedom (that is $\alpha$, $\beta$, $\gamma$, $\delta$, $c$). 
Therefore, the isospectrality will be
judged on the basis of very few eigenvalues. Our discretized operators 
correspond to $4\times 4$ matrices, so that we can only account on
four eigenvalues. This is the maximum number that allows us to find
continuously connected subfamilies of isospectral domains. In fact,
having more eigenvalues to handle can lead to an unsolvable problem, since
we do not have enough parameters to deform the domains. Thus, we
will look for a curve in the five-dimensional space, in such a way four
nonlinear equations (expressing the coincidence of the eigenvalues in relation
to an initial given quadrilateral) have to be simultaneously satisfied.

A first basic approach consists in implementing the classical finite difference method to  approximate 
the operator $L$ defined (\ref{defL})  and related to the 
eigenvalue problem (\ref{lapl3}), (\ref{laplBC3}). We recall that, when the quadrilateral
$\hat{Q}$ coincides with the square $Q$, one has $L=\Delta$, $\hat u=u$ and the corresponding
eigenvalue problem is given by (\ref{lapl1}),  (\ref{laplBC1}).

We take a uniform grid with step size $h$ in the unit square $Q\cup\partial Q$, both in the $x$ and $ y$ directions.
In particular, we divide the interval $[0,1]$  into three  equal subintervals of length $\frac13$.
The corresponding grid points are defined  by: 
\begin{equation}\label{gridDF}
G: \quad \quad ( x_{i}, y_{j})=(hi ,hj) = (\textstyle{\frac13} i, \textstyle{\frac13}j), \quad i,j=0,1,2,3. 
\end{equation}
Let $p_{1}=\left(\frac13, \frac13 \right)$, $p_{2}=\left(\frac23, \frac13 \right)$,
 $p_{3}=\left(\frac13, \frac23 \right)$, $p_{4}=\left(\frac23, \frac23 \right)$ be  the internal  points of the grid $G$.

In these circumstances using centered finite differences and  taking into account  the vanishing  Dirichlet boundary condition  
on the  boundary of $Q$ (see (\ref{laplBC3})),  
we have that the differential
operators $\frac{\partial^2}{\partial x^2}$,  $\frac{\partial^2}{\partial x \partial y}$, $\frac{\partial^2}{\partial y^2}$,
$\frac{\partial}{\partial x}$ and $\frac{\partial}{\partial y}$, 
can be,  respectively,   approximated with the discrete  operators
$D_{xx,G}$, 
$D_{xy,G}$, 
$D_{yy,G}$, 
$D_{x,G}$ and 
$D_{y,G}$
given by the following  $4\times 4$ matrices:

\begin{eqnarray} 
&&\hskip-0.8truecm D_{xx,G} =  
 9  \begin{pmatrix}
 -2   & \  1 \ &  0               &  0              \\
\  1 \ & -2 & 0              &  0              \\
 0                & 0              & -2 & \ 1 \  \\
0                & 0              & \ 1 \ & -2  \\
\end{pmatrix}, \qquad
D_{xy,G}=
\frac94  \begin{pmatrix}
0   & 0 & 0  & \ 1  \  \\
0   & 0 & \ -1 \  & 0   \\
0   & \ -1 \ & 0  & 0   \\\
\ 1  \ & 0 & 0  & 0   \\
\end{pmatrix}, \nonumber
\end{eqnarray}
\begin{eqnarray} \label{matrici}
&&\hskip-0.8truecm 
D_{yy,G}=
9  \begin{pmatrix}
-2   & 0 &\ 1   \            & 0              \\
0  & -2 & 0              & \ 1  \            \\
\ 1    \           & 0              & -2 & 0  \\
0                & \ 1    \         & 0&- 2  \\
\end{pmatrix}, \qquad
D_{x,G}=
\frac32  \begin{pmatrix}
0   & \ \ 1 \ \ & 0  & 0   \\
\ -1 \  & 0 & 0  & 0   \\
0   & 0 & 0  & \ \ 1 \ \   \\\
0   & 0 & \ -1  \ & 0   \\
\end{pmatrix},
\end{eqnarray}
\begin{eqnarray} 
&& 
D_{y,G}=
\frac32  \begin{pmatrix}
0   & 0& \ \ 1  \ \ & 0   \\
0   & 0 & 0  & \ \  1 \ \   \\
\ -1  \ & 0 & 0  & 0   \\\
0   & \ -1 \ & 0  & 0   \\
\end{pmatrix}.\nonumber
\end{eqnarray}
%
%
%
Finally, given  the functions $f_{i}$, $i=1,2,3,4,5,$   in (\ref{fi})  and defined  the five $4\times 4$ matrices:
\begin{eqnarray} \label{fiMatrici}
&&F_{i}=
  \begin{pmatrix}
f_{i}(p_{1})  & 0& 0  & 0   \\
0   & f_{i}(p_{2}) & 0  & 0   \\
0  & 0 & f_{i}(p_{3})  & 0   \\\
0   & 0 & 0  & f_{i}(p_{4})   \\
\end{pmatrix}, \qquad i=1,2,3,4,5, 
\end{eqnarray}
we have that the finite differences operator  $L_{fd}$ approximating  the differential 
operator  $L$  on the grid $G$ can be  written as follows:
\begin{eqnarray}\label{approxLsuG}
&& \hskip-1truecm L_{fd} =  F_{1} \cdot D_{xx,G} + F_{2} \cdot D_{xy,G} +F_{3}\cdot D_{yy,G} +F_{4} 
\cdot D_{x,G}+F_{5} \cdot D_{y,G},
\end{eqnarray}
where $\cdot$ denotes the usual matrix multiplication.
In this way the discrete  eigenvalue problem  associated to problem (\ref{lapl3}), (\ref{laplBC3}) can be written as:
\begin{equation}\label{eig3}
 L_{fd} \, \vec v =\lambda_{fd} \, \vec v  ,
\end{equation}
where the eigenvector $\vec v$ belongs to $\R^4$.

Better discretizations of $ \Delta$ and $L$ are  obtained by spectral collocation. In this simple circumstance one can
use polynomials of degree three in each variable. They must satisfy the vanishing constraint on the 
boundary of $Q$, so that one easily checks that the dimension of the approximating space is reduced
to four degrees of freedom. Thus, we collocate equation (\ref{lapl3}) at the four 
inner grid points in order to close the system.

 Indeed, we can also generalize the grid by introducing a new
parameter  $\kappa$,  with $0<\kappa<\frac12$,  and by  considering the following  set of points in the unit square $Q\cup\partial Q$:
\begin{equation}\label{gridDFn}
\bar{G}: \quad \quad (\bar{x}_{i}, \bar{y}_{j}), \quad i,j=0,1,2,3, 
\end{equation}
where 
\begin{equation}\label{gridDFpoint}
 (\bar{x}_{0},\bar{x}_{1},\bar{x}_{2},\bar{x}_{3})=(\bar{y}_{0},\bar{y}_{1},\bar{y}_{2},\bar{y}_{3})=(0,\kappa, 1-\kappa,1) .
\end{equation}
Of course, the internal nodes are:
 $\bar{p}_{1}=\left(\kappa, \kappa \right)$, $\bar{p}_{2}=\left(1-\kappa, \kappa \right)$,
 $\bar{p}_{3}=\left(\kappa, 1-\kappa \right)$, $\bar{p}_{4}=\left(1-\kappa, 1-\kappa \right)$.
They correspond to the classical uniform grid when $\kappa =1/3$. Another suitable choice is
to set $\kappa =\frac12-\frac{1}{2\sqrt5}$. In this way the values $2\kappa -1$ and $1-2\kappa$
are the zeros of the first derivative of the Legendre polynomial $P_3$.

Let us now consider a generic polynomial of degree 3 in each of the two variables $x$ and $y$:
\begin{equation}\label{pol}
{\cal P}(x,y)=a_{11} l_{1}(x) l_{2}(y) +  a_{12} l_{1}(x) l_{2}(y)+a_{21} l_{2}(x) l_{1}(y)+a_{22} l_{2}(x) l_{2}(y),
\end{equation}
for some coefficients $a_{ij} \in \R$, $i,j=1,2$.
In (\ref{pol})  $l_{1}$,   $l_{2}$   are  the one-dimensional Lagrange polynomials of degree 3  with respect 
to the nodes: $0,\kappa, 1-\kappa,1$. In particular $l_1$ and $l_2$ vanish at the endpoints, and we have:
\begin{equation}\label{Lag}
l_{1}(x)=\frac{x(x-1)(x+\kappa-1)}{\kappa (\kappa-1) (2\kappa-1)}, \quad \quad l_{2}(x)=-\frac{x(x-1)(x-\kappa)}{\kappa (\kappa-1) (2\kappa-1)}.
\end{equation}
By calculating explicitly the derivative of ${\cal P}$ and collocating at the grid points 
in  (\ref{gridDFn}), one obtains the  $4\times 4$ matrices: 
$D_{xx,\bar{G}}$, 
$D_{xy,\bar{G}}$, 
$D_{yy,\bar{G}}$,  
$D_{x,\bar{G}}$ and 
$D_{y,\bar{G}}$, 
approximating,  respectively,  the  differential operators 
 $\frac{\partial^2}{\partial x^2}$,  $\frac{\partial^2}{\partial x \partial y}$, $\frac{\partial^2}{\partial y^2}$,
$\frac{\partial}{\partial x}$ and $\frac{\partial}{\partial y}$.


Finally, as in (\ref{approxLsuG}), a spectral discretization $L_{sp}$ of the operator $L$ in the domain $Q$ takes the
form:
\begin{eqnarray}\label{approxLsuGbar}
&&\hskip-1truecm  L_{sp} =  F_{1} \cdot D_{xx,\bar{G}} + F_{2} \cdot D_{xy,\bar{G}} +F_{3} \cdot D_{yy,\bar{G}} +F_{4} \cdot
 D_{x,\bar{G}}+F_{5} \cdot D_{y,\bar{G}}, 
\end{eqnarray}
where the matrices   $F_{i}$, $ i=1,2,3,4,5,$ are defined in (\ref{fiMatrici}). This leads us
to the eigenvalue problem:
\begin{equation}\label{eig4}
 L_{sp} \, \vec v =\lambda_{sp} \, \vec v   , 
\end{equation}
where the eigenvector $\vec v$ belongs to $\R^4$.

We now discuss some practical cases. Suppose that $\hat{Q}=Q$, i.e., we are dealing with problem (\ref{lapl1}),  
(\ref{laplBC1}). From  (\ref{aut1Q}) we have that the first  four exact eigenvalues of minus the Laplace operator 
$(-\Delta )$ in the square are:
$(2\pi^2, 5\pi^2, 5\pi^2, 8\pi^2)\approx (19.74, 49.35, 49.35, 78.96)$. Regarding the
discretization of $-\Delta$ by finite differences  we find:
$(\lambda_{fd,1}, \lambda_{fd,2}, \lambda_{fd,3}, \lambda_{fd,4}) =(18,36,36,54)$. These values coincide
with those of the spectral approximation  for  $\kappa=\frac13$.
By using the spectral approximation with  $\kappa =\frac12-\frac{1}{2\sqrt5}$ we obtain instead $(\lambda_{sp,1}, 
\lambda_{sp,2}, \lambda_{sp,3}, \lambda_{sp,4}) =(20,40,40,60)$.

Taking into account the quadrilateral $\hat{Q}$  with 
vertices  $\hat{V}_{1}=(0,0)$, $\hat{V}_{2}=(1,0)$, $\hat{V}_{3}=(-0.2 ,1.1)$, $\hat{V}_{4}=(1.2 ,1.3)$ shown  in  
Fig.~\ref{fig2}
we find that the solution of  (\ref{eig3})  produces the eigenvalues 
$({\lambda}_{fd,1},{\lambda}_{fd,2}, {\lambda}_{fd,3}, {\lambda}_{fd,4}) =(12.54, 24.79, 25.43,  38.30)$.
Approximating $L$ with  (\ref{approxLsuGbar})  and  $\kappa=\frac13$  and solving  (\ref{eig4}) we find the following 
set of positive eigenvalues: $({\lambda}_{sp,1}, {\lambda}_{sp,2}, {\lambda}_{sp,3}, {\lambda}_{sp,4}) = 
(12.52, 24.63, 25.98,  38.05)$, while using the spectral approximation  (\ref{eig4})  with  $\kappa =\frac12-\frac{1}{2\sqrt5}$ 
we obtain $(\lambda_{sp,1}, \lambda_{sp,2}, $ $ \lambda_{sp,3}, \lambda_{sp,4}) =(13.92, 27.30,28.59,43.11)$.

\section{First attempts}\label{sec4}

\begin{figure}
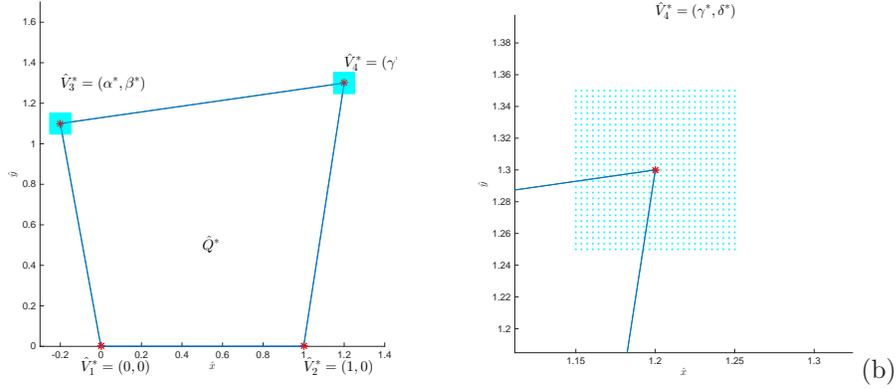

\centerline{\includegraphics[height=5cm]{fig3Griglie.eps}  (a)  \includegraphics[height=5cm]{fig3GriglieZoom.eps}  (b) }
  \caption{(a) The  quadrilateral $\hat{Q}^{*}$. Other possible quadrilaterals in the $h$-range
are obtained by selecting a point in each one of the two grids surrounding
the vertices on top. (b) Magnification of the area around the vertex $\hat{V}_{4}^{*}=(\gamma^{*} ,\delta^{*})$,
showing the right-hand side grid points.}
 \label{fig3}
  \end{figure}

Given the  family of quadrilateral  domains $\hat{Q} \subset \R^{2}$  of vertices
 $\hat{V}_{1}=(0,0)$, $\hat{V}_{2}=(1,0)$, $\hat{V}_{3}=(\alpha ,\beta)$, $\hat{V}_{4}=(\gamma ,\delta)$, 
 we begin our study  by fixing an initial  quadrilateral  $\hat{Q}^{*}$ of vertices $\hat{V}_{1}^{*}=\hat{V}_{1}=(0,0)$, 
$\hat{V}_{2}^{*}=\hat{V}_{2}=(1,0)$, $\hat{V}_{3}^{*}=(\alpha^{*} ,\beta^{*})$, $\hat{V}_{4}^{*}=(\gamma^{*} ,\delta^{*})$   
with $c=1$  (see  Fig.~\ref{fig3}).
From now on the superscript $*$
is used to point out that the domain $\hat{Q}^{*}$ has been fixed, and we would like to examine
what happens in its neighborhood in terms of isospectrality.

As shown in Sect.~\ref{sec3}, the quadrilateral $\hat{Q}^{*}$ is  mapped to the referring square $Q$ producing
the new operator $L$. By using the finite difference approximation  (\ref{approxLsuG}) or spectral collocation 
as in  (\ref{approxLsuGbar}),  the  discrete operator  can be either represented by
the  $4\times 4$ matrix $ L_{fd} $ or by $ L_{sp}$. Correspondingly,   problems  (\ref{eig3}) or  (\ref{eig4}) must 
be solved.

 From now on, we only work with $L_{sp}$ with $\kappa=\frac13$ and,  by solving (\ref{eig4}), we find four 
positive eigenvalues, that will be denoted by
${\lambda}_{sp,1}^{*}, {\lambda}_{sp,2}^{*}, {\lambda}_{sp,3}^{*}, {\lambda}_{sp,4}^{*}$. 
For example, as anticipated in the previous section,
 by taking the  quadrilateral $\hat{Q}^{*}$ of vertices 
 \begin{equation}\label{verticiQ*}
 \hat{V}_{1}^{*}=(0,0), \ \ \hat{V}_{2}^{*}=(1,0), \ \ \hat{V}_{3}^{*}=(-0.2 ,1.1), \ \ \hat{V}_{4}^{*}=(1.2 ,1.3),
  \end{equation}
%
shown in  Fig.~\ref{fig2} (see also Fig.~\ref{fig3}) and solving
 (\ref{eig4})   with $\kappa=\frac13$,    we obtain the following outcome:
\begin{equation}\label{autovreali}
({\lambda}_{sp,1}^{*}, {\lambda}_{sp,2}^{*}, {\lambda}_{sp,3}^{*}, {\lambda}_{sp,4}^{*}) =( 12.52, 24.63, 25.98,  38.05).
\end{equation}

Our goal is to see if there exist other quadrilaterals leading to the same
set of eigenvalues and if these can be connected by a curve.
The reason why we did not start our analysis with the unit square $Q$
and its isospectral companions will be clarified later on in Sect.~\ref{sec7}.

Since we do not have at the moment any theoretical result, a rough way to have an
idea of what happens it to check methodically a great number of 
quadrilaterals in the neighborhood of the initial one. To this purpose,
we construct two little squares with sides of length $l$  centered at the points $(\alpha^{*} ,\beta^{*})$ and
$(\gamma^{*} ,\delta^{*})$ as  shown  in Fig.~\ref{fig3}.  After defining appropriate grids of given size $h$
in these two squares, we try all the possible combinations. For each
couple of grid-points (one in the first square and one in the second square)  
we have a quadrilateral (recall that $\hat{V}_{1}^{*}=(0,0)$ and $\hat{V}_{2}^{*}=(1,0)$ have
been fixed). From this we deduce a $4\times 4$ matrix, and finally
four eigenvalues. We call the set of all quadrilaterals obtained
in this fashion: $h$-range.

We then select those configurations in the $h$-range displaying the same
eigenvalues given in (\ref{autovreali}), up to a prescribed error $\epsilon$. We can call
these special domains $\epsilon$-isospectral. The sizes of $h$
and $\epsilon$ have to be set up with the aim of finding reasonable
outcomes. In fact, if $\epsilon$ is too large we may end up with too many
$\epsilon$-isospectral domains; otherwise for $\epsilon$ too small
we could discover that the only acceptable domain is the starting one $\hat{Q}^{*}$.
Similar situations could also occur by selecting $h$ inappropriately.
The procedure is quite costly, especially for small $h$ and $\epsilon$.
This is the reason why in the next sections we look for something more 
convenient from the numerical viewpoint.

Unfortunately, the results of this analysis are not encouraging. Indeed,
it seems that there are no enough degrees of freedom to play with, and
this is the reason why in Sect.~\ref{sec2} we introduced the new parameter $c$. This is
an amplification (or reduction) factor that allows us to include
in the set of possible $\epsilon$-isospectral candidates other
quadrilaterals. These are obtainable through a suitable homothety centered
in $(0,0)$. 
Through this homothety  we associate
a generic point $\hat p$ of the plane $(\hat x, \hat y)$ to the point $\sqrt{c}  \ \hat p$. 
In this way, any given eigenvalue $\lambda_{sp}$ of $L_{sp}$ relative to a generic  domain $\Omega$ takes
the form $\lambda_{sp} /c$ in the new domain $\sqrt{c} \ \Omega$. By this trick, without too much effort, we can
enlarge the set of quadrilaterals at our disposition. 

 \begin{figure}
\centerline{\includegraphics[height=5cm]{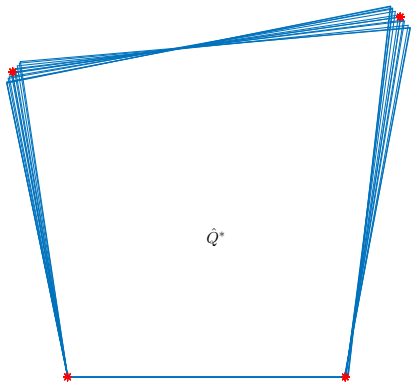} }
  \caption{Shape of some of the $\epsilon$-isospectral  domains related to the initial quadrilateral  
$\hat{Q}^{*}$ of Fig.~\ref{fig3} (a). The coordinates of their vertices have been previously 
multiplied by $\sqrt{c}$. In this fashion, for $c\not =1$, the second vertex on bottom
does not coincide with $(1,0)$.}
 \label{fig4}
  \end{figure}

Thus, we argue as follows. For a given domain in the $h$-range,
we will also accept eigenvalues that are proportional to those
in (\ref{autovreali}) through a multiplicative constant depending on $c$. 
To this end we set $ c={\lambda}_{sp,1}/{\lambda}_{sp,1}^{*}$. Afterwards
the 
domains are going to be selected according to the formula:
\begin{equation}\label{Errautov}
\left(\frac{  \sum_{k=1}^4  ({\lambda}_{sp,k}- c  \ {\lambda}_{sp,k}^{*})^{2}}
{  \sum_{k=1}^4 ({\lambda}_{sp,k}^{*})^{2}  }\right)^{1/2}
 \leq\epsilon .
\end{equation}
If the above inequality is satisfied for $c=1$, the corresponding
quadrilateral is directly $\epsilon$-isospectral to the starting one.
If $c\not =1$, then the actual $\epsilon$-isospectral quadrilateral is
obtained by multiplying the coordinates by the constant $\sqrt{c}$.

 For example, starting from the quadrilateral 
$\hat{Q}^{*}$  in (\ref{verticiQ*}) shown in   Fig.~\ref{fig3}  (which eigenvalues are listed  in (\ref{autovreali})),
we adopted the  above procedure based on the parameters:  $l=0.1$, $h=0.0036$ and $\epsilon=10^{-4}$. 
A family of  47  $\epsilon$-isospectral  domains  in the $h$-range was obtained. Some of these 
are displayed in Fig.~\ref{fig4}. Concerning the operator $L_{fd}$ in (\ref{approxLsuG}), preliminary results of this type
were found in \cite{codeluppi}.

We checked areas and perimeters of the   47  $\epsilon$-isospectral  domains. 
Within a tolerance of $10^{-3}$, 46 out of
47 quadrilaterals (including  $\hat{Q}^{*}$) share the same area, and 29 have the
same perimeter. Note that, in our discrete case, we cannot rely on a result
similar to that of Weyl for the continuous case. Nevertheless, the discovery that areas
are (almost) preserved, besides being in agreement with predictions, is an excellent tool to
decide a priori if a domain is appropriate. Indeed, before directly computing the
eigenvalues, one can filter those domains that, up to a certain accuracy, share
with the initial one the same area. This preliminary control saves a lot of computational
time.

\begin{figure}
\centerline{\includegraphics[height=4cm]{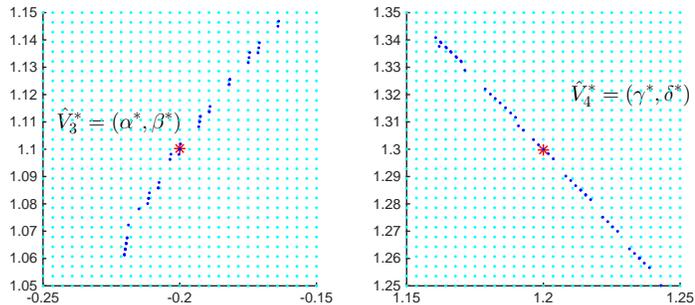} }
  \caption{Zoom of the little squares of Fig.~\ref{fig3} (a). The marked points are the
vertices of the  $\epsilon$-isospectral  domains starting from the quadrilateral 
$\hat{Q}^{*}$  in (\ref{verticiQ*}).}
 \label{fig5}
  \end{figure}

Note that the vertex $\hat{V}_{1}=(0,0)$ remains fixed, while
the vertex $\hat{V}_{2}$ shifts horizontally. This means that, except for the initial configuration
where the parameter $c$ is equal to one, the values of $c$ are in general different from one. 
The plots in Fig.~\ref{fig5} are the zooms of the two little squares     with sides of length $l=0.1$    
of Fig.   \ref{fig3} centered at 
the points $(\alpha^{*} ,\beta^{*})=(-0.2 ,1.1)$ and $(\gamma^{*} ,\delta^{*})=(1.2 ,1.3)$.
From these images we can conjecture the existence of a continuous path joining the
various  $\epsilon$-isospectral domains. Although only based on heuristic considerations,
the guess seems to be confirmed by further tests, where $h$ and $\epsilon$ are
conveniently taken smaller and smaller. We can be more precise in the coming sections, where
appropriate  strategies will be developed  to prove the existence of these curves and detect   them.

Similar results are obtained by varying the configuration of the
initial quadrilateral  $\hat{Q}^{*}$. We suggest however to set the initial parameters  
$ \alpha^{*} ,\beta^{*}, \gamma^{*} , \delta^{*}$ in order to stay away from some critical
situations that will be analyzed more in detail  in Sect.~\ref{sec7}.

\section{Implementation of the Implicit Function theorem}\label{sec5}

From the experiments of the previous section, our guess is that a curve
joining isospectral quadrilaterals actually exists. It is a one-parameter
function embedded in the 5 dimensional space spanned by the
parameters $\alpha$, $\beta$, $\gamma$, $\delta$, $c$. We can track it
by observing that it is implicitly characterized through four functional
equations. We can follow the curve by computing locally its tangent
vector and this can be done with the help of the Implicit Function theorem (see, for example, \cite{rudin}). 
Actually, in circumstances in which such a theorem is applicable, we
automatically have an existence result, at least at local level.

We recall that a generic quadrilateral 
$\hat{Q}$ has vertices of the form $\hat{V}_{1}=(0,0)$, 
$\hat{V}_{2}=(1,0)$, $\hat{V}_{3}=(\alpha ,\beta)$, $\hat{V}_{4}=(\gamma ,\delta)$.
Afterwords, for $i=1,2,3,4,$ let
 $\lambda_{sp,i}$ be  the eigenvalues in ascending order of the discrete problem  
(\ref{eig4}). Clearly, each one of these quantities depends on the parameters
$ \alpha ,\beta, \gamma , \delta$, i.e.:
\begin{eqnarray}\label{lambdaabcd}
&& \lambda_{sp,i}=\lambda_{sp,i} ( \alpha ,\beta, \gamma , \delta) , \quad  i=1,2,3,4.
\end{eqnarray}
Now, the eigenvalues (\ref{lambdaabcd}) can be seen as the roots
of the following characteristic polynomial:
%
\begin{eqnarray}\label{charpolyQ}
&& q(z)=z^{4} +\xi_{3}  z^{3} +\xi_{2}  z^{2} +\xi_{1}  z +\xi_{0}  ,
\end{eqnarray}
where  the coefficients   $\xi_{k}$, $k=0,1,2,3$, are computed from the entries of the matrix  $L_{sp}$ defined in (\ref{approxLsuGbar}). As a consequence these coefficients also depend on $ \alpha ,\beta, \gamma , \delta$, that is 
$\xi_{k}=\xi_{k}  ( \alpha ,\beta, \gamma , \delta)$,  $k=0,1,2,3$.

As done in the previous section, we can work around an initial domain $\hat{Q}^{*}$ 
of fixed vertices $\hat{V}_{1}^{*}=\hat{V}_{1}=(0,0)$, 
$\hat{V}_{2}^{*}=\hat{V}_{2}=(1,0)$, $\hat{V}_{3}^{*}=(\alpha^{*} ,\beta^{*})$, $\hat{V}_{4}^{*}=(\gamma^{*} ,\delta^{*})$   
with $c=1$. For instance, we  can start from the quadrilateral 
$\hat{Q}^{*}$  specified in (\ref{verticiQ*}) and  shown in   Fig.~\ref{fig3},  whose  corresponding eigenvalues are listed in  (\ref{autovreali}). Accordingly, $\lambda^{*}_{sp,i}$, $i=1,2,3,4,$ are
 the roots of the following characteristic polynomial:
\begin{eqnarray}\label{charpolyQ*}
&& q^{*}(z)=z^{4} +\xi_{3}^{*}  z^{3} +\xi_{2}^{*}  z^{2} +\xi_{1}^{*}  z +\xi_{0}^{*},
\end{eqnarray}
where now
the coefficients  $\xi_{k}^{*}$, $k=0,1,2,3$,  are given real numbers.
For example, starting  from the quadrilateral 
$\hat{Q}^{*}$  in (\ref{verticiQ*}),   we have:
\begin{eqnarray}
(\xi_{3}^{*}, \xi_{2}^{*}, \xi_{1}^{*},\xi_{0}^{*}) =( 101.18 ,    3675.65,     56468.45,     304819.78).
\end{eqnarray}

Thus, a certain quadrilateral 
$\hat{Q} \not =\hat{Q}^{*}$ (with vertices $\hat{V}_{1}=(0,0)$, 
$\hat{V}_{2}=(1,0)$, $\hat{V}_{3}=(\alpha ,\beta)$, $\hat{V}_{4}=(\gamma ,\delta)$) 
is isospectral to $\hat{Q}^{*}$ if and only if $q$ and $q^{*}$ have the same roots, i.e.,
one has $\xi_k=\xi_{k}^{*}$, $k=0,1,2,3$. We can weaken this condition by introducing the
parameter $c>0$. Indeed, we can also accept situations where the eigenvalues are proportionally related
as follows:
\begin{eqnarray}\label{charpolyQprimoequal0}
c=\frac{\lambda_{sp,i} }{\lambda_{sp,i}^{*}} , \quad i=1,2,3,4.
\end{eqnarray}
If $c=1$, $\hat{Q}$ turns out to be directly isospectral to $\hat{Q}^{*}$.
If $c\not =1$, the new domain  $\sqrt{c}  \ \hat{Q}$, obtained by homothety, is also isospectral to
$\hat{Q}^{*}$ (note that in this case  the second vertex $\hat{V}_{2}$  becomes $(\sqrt{c},0)$). We can now translate condition 
(\ref{charpolyQprimoequal0}) in terms of polynomial coefficients, obtaining:
\begin{eqnarray}\label{polcof}
c^{4-k}=\frac{\xi_{k} }{\xi_{k}^{*}} , \quad k=1,2,3,4.
\end{eqnarray}

In the end, we propose to introduce the four  functions ${\cal F}_{k}={\cal F}_{k}(\alpha, \beta, \gamma, \delta, c)$, $ k=0,1,2,3,$  of  the five variables $\alpha, \beta, \gamma, \delta, c$,  defined as follows:
\begin{eqnarray}\label{Funzionik}
{\cal F}_{k}(\alpha, \beta, \gamma, \delta, c)= \xi_{k}  -c^{4-k} \ \xi_{k}^{*} , \qquad k=0,1,2,3,
\end{eqnarray}
and look for values such that: 
\begin{eqnarray}\label{Funzionikzero}
{\cal F}_{k}=0, \qquad  k=0,1,2,3.
\end{eqnarray}
 If we are lucky, there
is a local curve $\Phi (t)\in \R^5$, $t\in [-T, T]$, $T>0$,  described by the functions 
$\alpha (t)$, $\beta (t)$, $\gamma (t)$, $\delta (t)$, $c(t)$, passing through the point 
${\bf P}^*$ of coordinates $\alpha^{*}$, 
$\beta^{*}$, $\gamma^{*}$, $\delta^{*}$, $c=1$ (i.e., the parameters identifying  the initial 
quadrilateral $\hat{Q}^{*}$ as in   (\ref{verticiQ*})) and connecting isospectral domains.

In order to  follow such a curve  $\Phi $, we need to find its local tangent vector. We can 
express the various parameters in function of one of them.
We fix for example $\beta (t)=\beta^{*}+t$, with $t\in  [-T, T]$, so that $\beta (0)=\beta^{*}$. 
We then differentiate with
respect to $t$ the four equations in (\ref{Funzionikzero}), arriving at the system:
%
%
\begin{eqnarray}
&&
  \begin{pmatrix}
  \ \ \displaystyle\frac{\partial {\cal F}_0}{\partial \alpha} \ \ &  \ \ \displaystyle \frac{\partial {\cal F}_0}{\partial \gamma}  \ \    
		& \ \  \displaystyle\frac{\partial {\cal F}_0}{\partial \delta} \ \  &  \ \ \displaystyle\frac{\partial {\cal F}_0}{\partial c}   \ \        \\[4mm]
 \displaystyle\frac{\partial {\cal F}_1}{\partial \alpha}  &  \displaystyle \frac{\partial {\cal F}_1}{\partial \gamma}        &   \displaystyle\frac{\partial {\cal F}_1}{\partial \delta} &  \displaystyle\frac{\partial {\cal F}_1}{\partial c}           \\[4mm]
  \displaystyle\frac{\partial {\cal F}_2}{\partial \alpha}  &  \displaystyle \frac{\partial {\cal F}_2}{\partial \gamma}        &   \displaystyle\frac{\partial {\cal F}_2}{\partial \delta} &  \displaystyle\frac{\partial {\cal F}_2}{\partial c}           \\[4mm]
   \displaystyle\frac{\partial {\cal F}_3}{\partial \alpha}  &  \displaystyle \frac{\partial {\cal F}_3}{\partial \gamma}        &   \displaystyle\frac{\partial {\cal F}_3}{\partial \delta} &  \displaystyle\frac{\partial {\cal F}_3}{\partial c}           \\
\end{pmatrix}
  \begin{pmatrix}
  \displaystyle\frac{d \alpha}{d t}     \\[4mm]
  \displaystyle\frac{d \gamma}{d t}     \\[4mm]
  \displaystyle\frac{d \delta}{d t}     \\[4mm]
  \displaystyle\frac{d c}{d t}     \\
\end{pmatrix}
=
  \begin{pmatrix}
-  \displaystyle\frac{\partial  {\cal F}_0}{\partial \beta}     \\[4mm]
-  \displaystyle\frac{\partial  {\cal F}_1}{\partial \beta}     \\[4mm]
-  \displaystyle\frac{\partial  {\cal F}_2}{\partial \beta}     \\[4mm]
-  \displaystyle\frac{\partial  {\cal F}_3}{\partial \beta}     \\
\end{pmatrix} ,\label{problemFk}
\end{eqnarray}
in the unknowns $\alpha^{\prime} (t)=\frac{d \alpha}{d t}$,  
$\gamma^{\prime} (t)=\frac{d \gamma}{d t}$,  $\delta^{\prime} (t)=
\frac{d \delta}{d t}$ and $c^{\prime} (t)=\frac{d c}{d t}$.

At this point, it is important to observe that the functions  ${\cal F}_{k}$, $k=0,1,2,3$, in   (\ref{Funzionik})
and their derivatives are explicitly known, although their expressions may result rather
complicated. In fact, starting from $\alpha$, $\beta$, $\gamma $, $\delta$, one
can build the coefficients of the mapping into the referring square.  Successively,
always in function of these parameters,  one writes the matrix $L_{sp}$ defined in (\ref{approxLsuGbar}). Finally,
a closed form is also known for the coefficients of the characteristic polynomial 
(this is not true instead for the eigenvalues). Of course, from the practical viewpoint, these computations
can only be carried out with the help of a software running with symbolic manipulation.

Going through this calculation, we find that the determinant of the matrix in
(\ref{problemFk}) is different from zero at ${\bf P}^*$. Therefore, one is able to theoretically detect a 
value $T>0$, such that a curve  $\Phi (t)\in \R^5$, $t\in  [-T, T]$,  described by the parameters 
$\alpha (t)$, $\beta (t)=\beta^{*}+t$, $\gamma (t)$, $\delta (t)$, $c(t)$, actually exists
in the interval $t\in  [-T, T]$. For $t=0$, we have $\Phi (0)=(\alpha^*,\beta^*, \gamma^*,
\delta^*,1)={\bf P}^*$.

The following  first-order explicit iteration method can be adopted in order to  
follow the curve, at least locally.
Recall that, given $T>0$, we have  chosen   $\beta (t)=\beta^{*}+t$,  $t\in  [-T, T]$.    
For simplicity let us fix our attention on positive value of  the variable $t$ (indeed
  similar arguments may also be applied when  $t$ is negative)
and, given an integer $M>0$,  let us discretize the interval $[\beta^{*},\beta^{*}+T]$ with a uniform grid
$t_{m}$, $m=0,\ldots, M$, with step size $\delta t=\frac{T}{M}$, that is:
\begin{equation}\label{grigliaBeta}
t_{m}=\beta^{*} +m \delta t=\beta^{*} +m  \frac{T}{M}, \quad          m=0,\ldots, M.
\end{equation}
Our curve $\Phi$ is going to be approximated by the quantities $\Phi_m \approx \Phi (t_m)$,
$m=0,\ldots, M$. In particular $\Phi_0 =\Phi (0) ={\bf P}^*$. At each step, the exact
derivatives $\alpha^{\prime} (t_m)$,  $\gamma^{\prime} (t_m)$,  $\delta^{\prime} (t_m)$
 and $c^{\prime} (t_m)$, $m=0,\ldots, M$, are computed by solving (\ref{problemFk}) (note that $\beta^\prime
(t_m)=1$, $m=0,\ldots, M$). These exact
derivatives are organized in the correcting vector $\Psi_m$, $m=0,\ldots, M$. The algorithm is
then based on the following  iteration:
\begin{equation}\label{iter}
\Phi_{m+1}=\Phi_m +\delta t \Psi_m \qquad  m=0,\ldots,M-1.
\end{equation}
A similar procedure can be used  to approximate the curve  $\Phi $ in the interval $ [-T, 0]$.

For example, we discuss the case of the quadrilateral 
$\hat{Q}^{*}$  in (\ref{verticiQ*}) shown in   Fig.~\ref{fig3} and  we take  $\beta (t)=\beta^{*}+t=1.1+t$,  $t\in  [-T, T]= [-0.06,0.06]$.  
We apply  (\ref{iter}) with $M=100$, both for positive and negative values of $t$. 
Fig.~\ref{fig6} and Fig.~\ref{fig7} show, respectively,  the projections  of $\Phi_{m}$, $m=0,\ldots, M$, in the planes 
$(\alpha ,\beta )$,  $(\gamma ,\delta )$ and the graph $(t_m,c(t_m))$,   $m=0,\ldots, M$.
These results confirm what was predicted in the previous section. 

As expected, the shapes of the projections  of the  curve turn out to be exactly the same if, instead of the parameter $\beta$, another
parameter is assumed to be explicit. 

 \begin{figure}
\centerline{\includegraphics[height=4cm]{fig6.eps} }
  \caption{Projections of the curve $\Phi (t)\in \R^5$, $t\in   [-T, T]= [-0.06,0.06]$,  in the
	planes $(\alpha,\beta )$ and $(\gamma, \delta )$. These have been found with the 
help of the Implicit Function theorem and  the iteration algorithm   (\ref{iter}).}
 \label{fig6}
  \end{figure}
 \begin{figure}
 \centerline{\includegraphics[height=4.3cm]{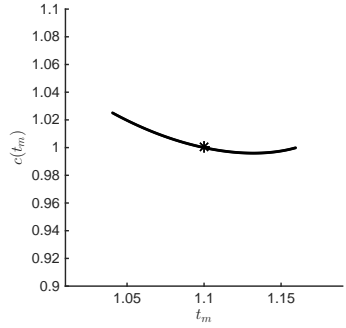} }
  \caption{Plot of the points $(t_m,c(t_m))$,   $m=0,\ldots, M$,  corresponding
 to the curves of Fig.~\ref{fig6}. We recall that $\sqrt{c}$ is an amplification  (or reduction) factor,
used to adjust the size of the domains, employed as an additional parameter to ensure
isospectrality. }
 \label{fig7}
  \end{figure}
  %

\section{Other approaches}\label{sec6}

The use of symbolic manipulation, in order to compute the exact tangent vector to
the curve $\Phi$ at a given point, is certainly expensive, ill-conditioned and
difficult to be generalized. It has been however an important theoretical tool
to establish the local existence of the curve.
In this section we propose an alternative numerical procedure, far more cheaper but
with similar performances.


In line with the arguments invoked  in Sect.~\ref{sec5} and using the same notation we 
have that,  if  the curve  $\Phi$ has to connect isospectral domains,   the eigenvalues  
$  \lambda_{sp,i}/c$, $ i=1,2,3,4$,  in  (\ref{charpolyQprimoequal0})   must remain the same
as $t$ varies.  In particular we have: $\lambda_{sp,i}^{*}= 
\lambda_{sp,i}/c$, $ i=1,2,3,4$.
That  is, having in mind  (\ref{lambdaabcd}),  the four functions  ${\cal G}_{i}={\cal G}_{i}(\alpha, \beta, 
\gamma, \delta, c)$, $  i=1,2,3,4,$  of  the five variables $\alpha, \beta, \gamma, \delta, c$, defined as follows:
\begin{eqnarray}\label{FunzionikG}
{\cal G}_{i}  (\alpha, \beta, \gamma, \delta, c)= \frac{\lambda_{sp,i} }{c}   , \quad  i=1,2,3,4, 
\end{eqnarray}
must be constant, i.e., their derivatives with respect to $t$ are zero.
By differentiating (\ref{FunzionikG}) with respect to $t$  and  by expressing the 
various parameters in function 
of one of them,  we find the approximated local tangent vector to the curve  $\Phi $.

As in Sect.~\ref{sec5}   let us  fix for example $\beta (t)=\beta^{*}+t$, with $t\in  [-T, T]$, so that $\beta (0)=\beta^{*}$. 
The new linear system takes the form:
%
%
\begin{eqnarray}
&&
  \begin{pmatrix}
\ c \   \displaystyle\frac{\partial  \lambda_{sp,1}}{\partial \alpha} \  & \  c \   \displaystyle\frac{\partial  \lambda_{sp,1}}{\partial \gamma} \ &
\ c \   \displaystyle\frac{\partial  \lambda_{sp,1}}{\partial \delta} \ & \ \ -  \lambda_{sp,1} \ \\[4mm]
c \   \displaystyle\frac{\partial  \lambda_{sp,2}}{\partial \alpha}  &  c \   \displaystyle\frac{\partial  \lambda_{sp,2}}{\partial \gamma} &
c \   \displaystyle\frac{\partial  \lambda_{sp,2}}{\partial \delta}  & -  \lambda_{sp,2} \\[4mm]
c \   \displaystyle\frac{\partial  \lambda_{sp,3}}{\partial \alpha}  &  c \   \displaystyle\frac{\partial  \lambda_{sp,3}}{\partial \gamma} &
c \   \displaystyle\frac{\partial  \lambda_{sp,3}}{\partial \delta}  & -  \lambda_{sp,3} \\[4mm]
c \   \displaystyle\frac{\partial  \lambda_{sp,4}}{\partial \alpha}  &  c \   \displaystyle\frac{\partial  \lambda_{sp,4}}{\partial \gamma} &
c \   \displaystyle\frac{\partial  \lambda_{sp,4}}{\partial \delta}  & -  \lambda_{sp,4} \\
\end{pmatrix}
%
  \begin{pmatrix}
  \displaystyle\frac{d \alpha}{d t}     \\[4mm]
  \displaystyle\frac{d \gamma}{d t}     \\[4mm]
  \displaystyle\frac{d \delta}{d t}     \\[4mm]
  \displaystyle\frac{d c}{d t}     \\
\end{pmatrix}
=
  \begin{pmatrix}
- c \   \displaystyle\frac{\partial  \lambda_{sp,1}}{\partial \beta}      \\[4mm]
- c \   \displaystyle\frac{\partial  \lambda_{sp,2}}{\partial \beta}      \\[4mm]
- c \   \displaystyle\frac{\partial  \lambda_{sp,2}}{\partial \beta}      \\[4mm]
- c \   \displaystyle\frac{\partial  \lambda_{sp,2}}{\partial \beta}      \\
\end{pmatrix} ,\label{problemGk}
\end{eqnarray}
in the unknowns $\alpha^{\prime} (t)=\frac{d \gamma}{d t}$,  
$\gamma^{\prime} (t)=\frac{d \gamma}{d t}$,  $\delta^{\prime} (t)=
\frac{d \delta}{d t}$ and $c^{\prime} (t)=\frac{d c}{d t}$.

Since now the explicit expression of the partial derivatives in (\ref{problemGk}) is not available,
the entries of the matrix are approximated by classical finite differences, that is,  for example, given an 
increment $d\alpha$, we can compute $\partial  \lambda_{sp,i} /\partial \alpha$, $i=1,2,3,4$, as follows:
\begin{eqnarray}\label{derivAlpha}
\displaystyle\frac{\partial  \lambda_{sp,i}}{\partial \alpha} \approx
\frac{\lambda_{sp,i}(\alpha+d\alpha, \beta, \gamma, \delta, c)-\lambda_{sp,i}(\alpha, \beta, \gamma, \delta, c)}{d\alpha}  , \quad  i=1,2,3,4.
\end{eqnarray}
Similar  formulae can be used to approximate the other entries of the matrix in (\ref{problemGk}).

Finally,  we use the uniform grid (\ref{grigliaBeta}) when $t$ is positive (or its opposite when $t$ is negative).
At step $m$, we evaluate the approximations of 
  $\alpha^{\prime} (t_m)$,  $\gamma^{\prime} (t_m)$,  $\delta^{\prime} (t_m)$
 and $c^{\prime} (t_m)$, $m=0,\ldots, M$  (note that $\beta^\prime
(t_m)=1$, $m=0,\ldots, M$), by solving a system similar to that in (\ref{problemGk}), where partial
derivatives have been replaced with incremental ratios. We put the so obtained
values in the current vector $ \bar{\Psi}_m$, $m=0,\ldots, M$. Thus, a rough approximation of  $\Phi (t_m)$  
is given  by the quantities $\bar{\Phi}_m $, $m=0,\ldots, M$, recovered by recursion
from the following explicit iteration method: 
\begin{equation}\label{iter2}
\bar{\Phi}_{k+1}=\bar{\Phi}_k +\delta t \  \bar{\Psi}_k \qquad  k=0,\ldots,M-1,
\end{equation}
starting from   $\bar{\Phi}_0 =\Phi (0) ={\bf P}^*$.

\begin{figure}
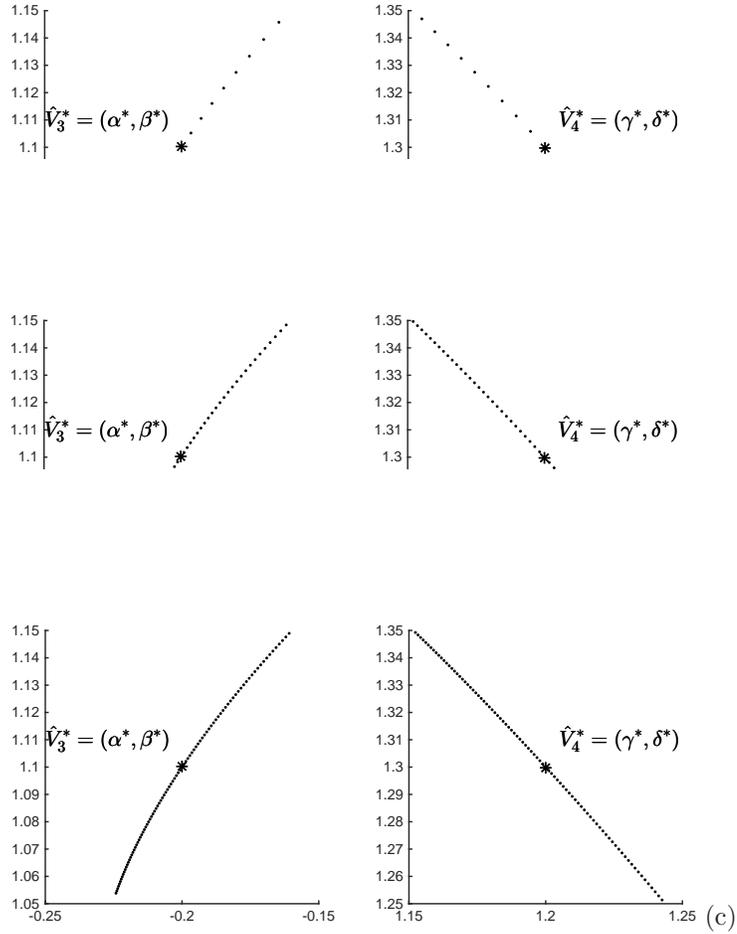

\centerline{\includegraphics[height=4cm]{fig8M10.eps} (a) }
\centerline{\includegraphics[height=4cm]{fig8M30.eps} (b) }
\centerline{\includegraphics[height=4cm]{fig8M50.eps} (c) }
  \caption{Projections of the curve $\Phi (t)\in \R^5$, $t\in   [-T, T]= [-0.06,0.06]$,  found with the 
help of the iteration algorithm   (\ref{iter2}) with $M=10$ (case (a)), $M=30$ (case  (b)), $M=50$ (case  (c)).
These graphs have to be compared with those of Fig.~\ref{fig6} }
 \label{fig8}
  \end{figure}

Given the quadrilateral  $\hat{Q}^{*}$  corresponding to (\ref{verticiQ*}) and shown in   Fig.~\ref{fig3},  
by taking  $\beta (t)=\beta^{*}+t=1.1+t$,  $t\in  [-T, T]= [-0.06,0.06]$,  
we apply  the iteration method (\ref{iter2}), both for positive and negative values of $t$.  Different 
values of the  step size $\delta t=\frac{T}{M}$ are used. In particular we choose $M=10, 30, 50$.
In Fig.~\ref{fig8}    we show  the projections of  $\bar{\Phi}_{m}$,   $m=0,\ldots, M$,  in 
the planes $(\alpha ,\beta )$,  $(\gamma ,\delta )$.
 As $\delta t$ diminishes, the graphs of Fig.~\ref{fig8} approach those of  Fig.~\ref{fig6}.
 
Also in this case the shape of the projections of the curve turn out to be qualitatively the same if
another parameter is assumed to be explicit (instead of the parameter $\beta$). 

This convergence behavior is a good surprise. Indeed, such a property should not be given for granted. In fact, 
through  (\ref{derivAlpha}) and similar formulae,  we replaced partial derivatives by a first-order approximation,
and we introduced this correction in the first-order algorithm (\ref{iter}). This double
discretization does not necessarily bring to a convergent scheme.

A sort of convergence analysis can be carried out by examining the history of the quadrilateral areas
in comparison to the area of  $\hat{Q}^{*}$, given by $\mu_2 (\hat{Q}^{*})=1.44$. The results of
some tests are reported in Table \ref{tavola}. For a fixed $\delta t$ the error
grows linearly with the distance from the quadrilateral  $\hat{Q}^{*}$.
Nevertheless, there is no substantial decay of the error by diminishing $\delta t$.
It has also to be noted however that, in the discrete case, we do not have any theoretical
result ensuring that areas of isospectral domains must be preserved. 

We finally observe
that the results obtained so far do not change significantly if other values of the
parameter $\kappa$ introduced in (\ref{gridDFpoint}) are taken into account (we recall that
in the experiments reported in this paper we used $\kappa =\frac13$).

\begin{table}
\caption{Relative errors of the areas in relation to the initial
position.}
\label{tavola}      
\begin{tabular}{llll}
\hline\noalign{\smallskip}
$\beta^{*}+t$ &  $M=5$ and $\delta t$ =0.012 &  $M=10$  and $\delta t $=0.006 &  $M=20$ and
  $\delta t$ =0.003\\
 \hline\noalign{\smallskip} 
1.040 &    11.70 e-04&     10.30 e-04   &   9.65 e-04\\
1.052 &    9.91 e-04&  8.69   e-04   &  8.12  e-04\\
1.064 &     7.83 e-04&   6.84  e-04   &  6.37  e-04\\
1.076  &    5.48 e-04&   4.76  e-04   &  4.43  e-04\\
1.088  &   2.87 e-04&    2.48 e-04   &   2.30  e-04\\
1.1 &    0&    0   &    0\\
1.112 &    1.53 e-04&   1.89  e-04   &   2.07 e-04\\
1.124 &    3.09 e-04&   3.85  e-04   &  4.25  e-04\\
1.136 &    4.67 e-04&    5.89 e-04   &   6.54 e-04\\
1.148  &    6.26 e-04&    8.00 e-04   &   8.92 e-04\\
 1.160 &   7.86 e-04&     10.19 e-04  &    11.42 e-04\\
\noalign{\smallskip}\hline
\end{tabular}
\vspace{.4cm}
\end{table}
%

\section{The case of the square}\label{sec7}
%

Here, we discuss about the case of the  unit square $Q$, i.e. the quadrilateral  of vertices  $V_{1}=(0,0)$, $V_{2}=(1,0)$, $V_{3}=(0,1)$, $V_{4}=(1,1)$,  and the case of  other domains  $\hat{Q}$ of 
vertices  $\hat{V}_{1}=(0,0)$, $\hat{V}_{2}=(1,0)$, $\hat{V}_{3}=(\alpha ,\beta)$, $\hat{V}_{4}=(\gamma ,\delta)$  that
are symmetric with respect to the straight line $y=x$, that is domains $\hat{Q}$ with  $\alpha =0$, $\beta
=1$ and $\gamma =\delta$.  Of course, when $\alpha =0$, $\beta
=1$, $\gamma =\delta=1$  we have   $\hat{Q}=Q$.

\begin{figure}
\centerline{\includegraphics[height=4cm]{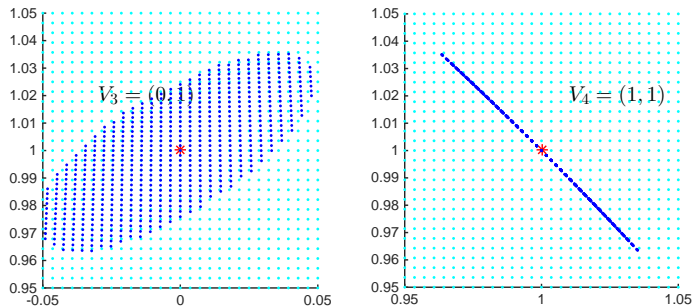} }
  \caption{Zoom of the vertices of the  $\epsilon$-isospectral  domains when the initial quadrilateral
is the unit square $Q$.}
 \label{fig9}
  \end{figure}

In these situations using the approach proposed in Sect.~\ref{sec5}, one can check that the determinant of the matrix
in (\ref{problemFk}) is always  zero.
We remind that  in Sect.~\ref{sec5}   we have chosen  $\beta (t) =1+t$, $t\in  [-T, T]$, but, indeed,  the determinant remains zero
independently on the explicited  parameter. 
Moreover, for the square  $Q$, the rank of the
matrix in (\ref{problemFk})  is just one, i.e. all the lines (or rows) of the matrix are linearly dependent.
Definitely, we are not in the situation that allows us to use the
Implicit Function theorem. Nevertheless, experimentally, one can find
in the neighborhood of the unit square $Q$ many other quadrilaterals isospectral to
it. 
An analysis similar to that of Sect.~\ref{sec4}  
reveals configurations as the ones shown in Fig.~\ref{fig9}. 
In this case we choose:  $l=0.1$, $h=0.0036$ and $\epsilon=5\times 10^{-4}$, and, as in  Sect.~\ref{sec4}, we 
considered $c={\lambda}_{sp,1}/{\lambda}_{sp,1}^{*}$.
In the neighborhood
of the vertex  $V_{3}=(0,1)$ on the left of Fig.~\ref{fig9} the distribution of selected points does not reveal specific
patterns, while for the vertex $V_{4}=(1 ,1)$ on the 
right of Fig.~\ref{fig9} there is a superimposition of several curves, as it will be evident from the
discussion that follows.

  %

%
 \begin{figure}
\centerline{\includegraphics[height=4.cm]{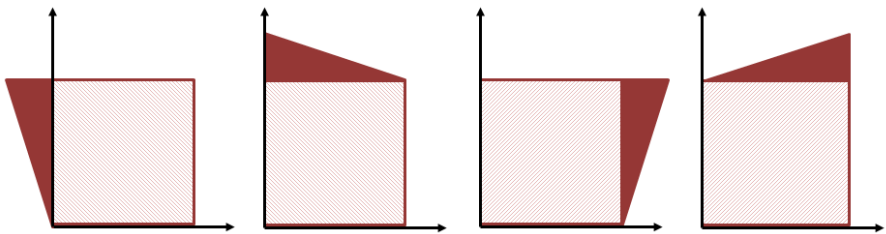} }
  \caption{Deformations of the unit square  $Q$ obtained by moving independently
each single parameter. The resulting quadrilaterals are isometric.}
 \label{fig10}
  \end{figure}

A heuristic explanation of the non applicability of the Implicit Function
theorem is inspired by Fig.~\ref{fig10}, where the four parameters
$-\alpha$,
$\beta$, $\gamma$, $\delta$, are varied independently (one has to imagine
that these deformations are infinitesimal). One gets four configurations
corresponding to quadrilaterals that are reciprocally isometric (by rotations
or symmetries), and therefore isospectral. This somehow explains why the
Jacobian matrix in  (\ref{problemFk}) has rank equal to one. Analogous
conclusions hold for other initial domains presenting symmetries.

 \begin{figure}
\centerline{\includegraphics[height=4cm]{fig11.eps} }
 \caption{Zoom of the vertices of the  $\epsilon$-isospectral  domains starting from the unit square $Q$ and  imposing $\alpha=0$ in the $h$-range.}
 \label{fig11}
%
\centerline{\includegraphics[height=4cm]{fig12.eps} }
 \caption{Zoom of the vertices of the  $\epsilon$-isospectral  domains starting from the unit square $Q$  and  imposing $\beta=1+\alpha$ in the $h$-range.}
 \label{fig12}
  \end{figure}
%
%
 \begin{figure}
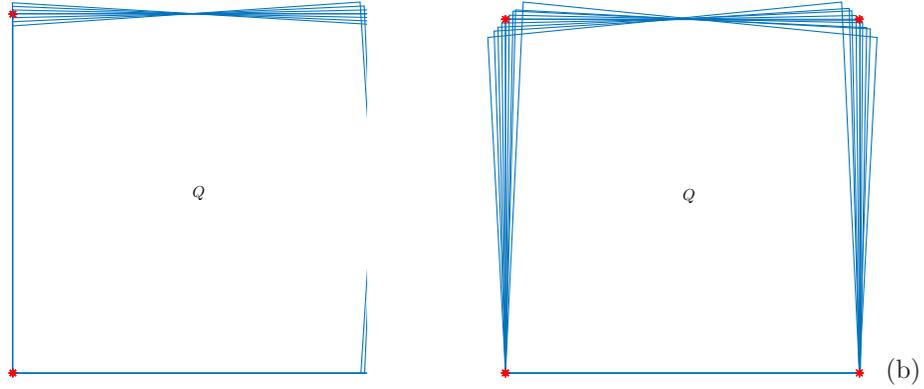

\centerline{\includegraphics[height=5cm]{fig13a.eps} (a)\hskip.8truecm\includegraphics[height=5cm]{fig13b.eps} (b)}
  \caption{Shapes of the $\epsilon$-isospectral  domains obtained by departing from the initial unit square  
$Q$, and  imposing $\alpha=0$ (case   (a)) or $\beta=1+\alpha$ (case (b)). }
 \label{fig13}
  \end{figure}

The above considerations do not prevent however the existence of curves connecting isospectral domains.
Actually, in correspondence of the parameters associated with the unit square  $Q$
we have a biforcation point. A deeper study shows that different curves of isospectral quadrilaterals are
obtained. For instance, we can have those of the family associated either with the points
of Fig.~\ref{fig11} or with the points of Fig.~\ref{fig12}. In  particular,  Fig.~\ref{fig13} shows the
quadrilaterals corresponding to the two ways of deforming the square  $Q$. These two branches
of isospectral domains departing from the unit square are very neat. Nevertheless, their
identification is not easy if one only examines Fig.~\ref{fig9}. Of course, when we started
our analysis we (erroneously) thought that the case of the square was the easiest one; only later we realized
that this was far from being true.

\begin{figure}
\centerline{\includegraphics[height=5cm]{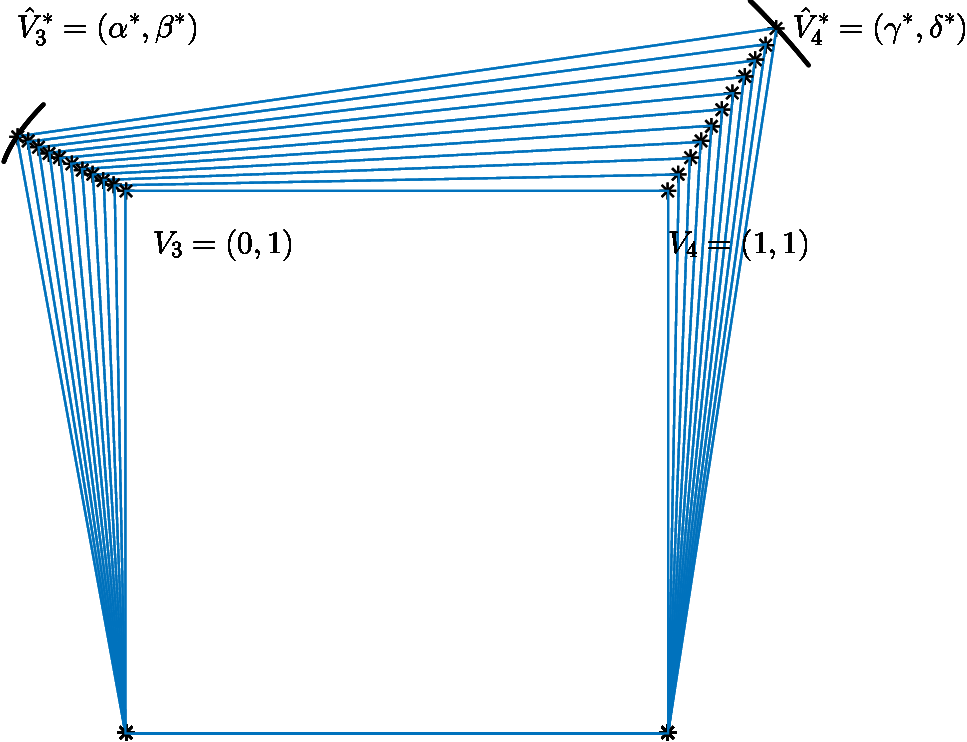} }
  \caption{ The quadrilateral  $\hat{Q}^{*}$  in (\ref{verticiQ*}) is transformed into the unit square $Q$ in 10 steps.}
	\label{fig14}
\vskip1truecm
\centerline{\includegraphics[height=4cm]{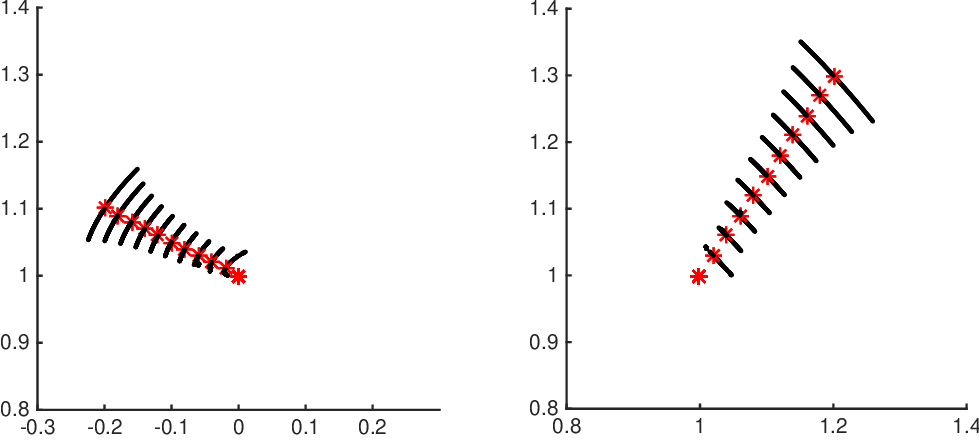} }
 \caption{Projections of the local curves joining the 
	corresponding isospectral domains for each of the quadrilaterals of Fig.~\ref{fig14}.
Note that for $j=9$, on the left-hand side the corresponding curve shows  a kind of turning point.
It is not easy to guess how the bifurcation point (emerging when $j=S=10$) is approached.}
 \label{fig15}
  \end{figure}

We conclude this section with a few more experiments. We would like to figure out what
happens to the curve connecting isospectral domains when the initial quadrilateral $ \hat{Q}^{*}$
is modified. For example, we can start from the quadrilateral $ \hat{Q}^{*}$  associated with the 
vertices in (\ref{verticiQ*}). The corresponding isospectral family forms a curve characterized 
by the projections shown in Fig.~\ref{fig6} (also reported in  Fig.~\ref{fig14}). Now, we deform 
such initial domain $ \hat{Q}^{*}$ by slowly approaching the unit square $Q$. To this purpose, we fix an integer $S>1$.
For  $j=0,\ldots, S$, the vertices of the
transitory quadrilaterals $\hat Q_j$ are chosen according to the law:  
\begin{eqnarray} \label{legge}
&& \hat{V}_{1,j}=V_{1}=(0,0), \quad \hat{V}_{2,j}=V_{2}=(1,0), \nonumber\\
&&  \hat{V}_{3,j}=(1-s_j)\hat V_{3}^{*} +s_j{V}_{3}, \quad \hat{V}_{4,j}=(1-s_j)\hat V_{4}^{*} +s_j{V}_{4}, \nonumber\\
&&\hskip6truecm \ s_j=j/S, \, \, j=0,\ldots, S.
\end{eqnarray}
In this way, for $j=0$ in (\ref{legge}) we have the initial quadrilateral $ \hat{Q}^{*}$, i.e.  
$\hat Q_0=\hat{Q}^{*}$,   while 
for $j=S$  we obtain the unit square $Q$, i.e.  $\hat Q_S=Q$. 

Since we know
that, in correspondence of the unit square  $Q$, the determinant of the Jacobian in (\ref{problemFk})
is zero, the use of the Implicit Function theorem becomes more and more restrictive as $j$ increases.
For this reason we need to limit the range $[-T_j, T_j]$ of definition of the curve.
Given $T_0 >0$ we propose to control this range  by defining the extremes 
$ T_j$,  $j=0,\ldots, S-1$,  of the interval according to the rule: 
\begin{equation}\label{legget}
T_j=\frac{T_0}{1+2*s_j},  \quad \quad  j=0,\ldots, S-1.
\end{equation}
We excluded the value  $j=S$  in  (\ref{legget})  because in this case  the procedure is not applicable.

Starting from the quadrilateral 
$\hat{Q}^{*}$  in (\ref{verticiQ*}), taking  $S=10$ in (\ref{legge}) and $T_0=0.06$ in  (\ref{legget}), we obtain  the various quadrilaterals 
displayed in Fig.~\ref{fig14}. For  $j=0,\ldots, S-1$, we locally compute the curve of the
isospectral domains related to $\hat Q_j$, indifferently with the  approach proposed in  Sect.~\ref{sec6} or Sect.~\ref{sec7}. 
The projections of these curves in the planes
$(\alpha ,\beta )$ and $(\gamma ,\delta )$ are shown in Fig.~\ref{fig15}. 

\begin{figure}
\centerline{\includegraphics[height=5cm]{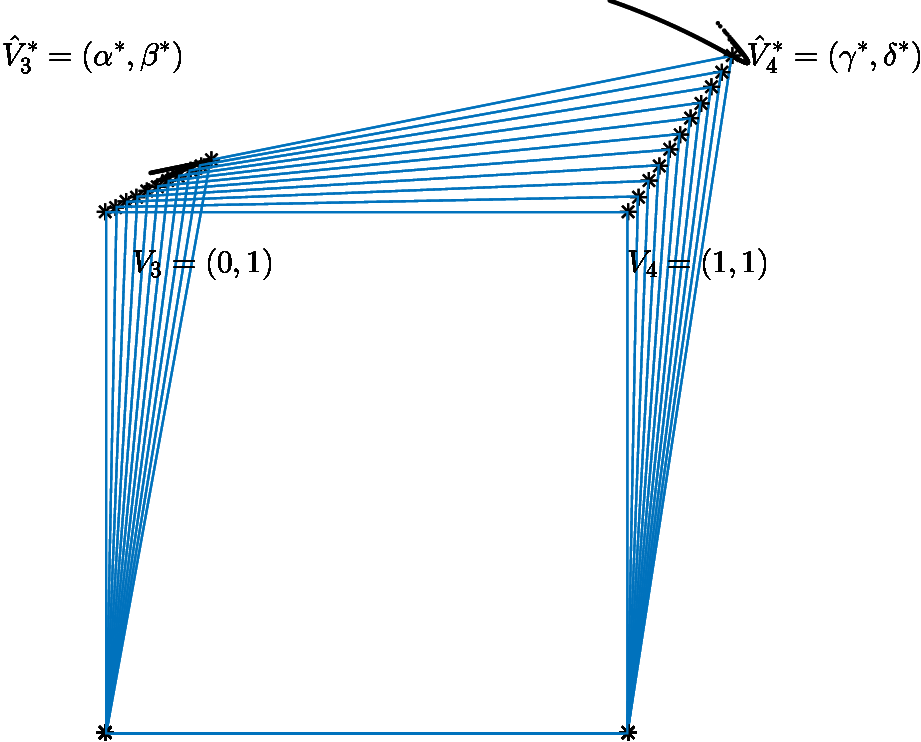} }
  \caption{The quadrilateral  $\hat{Q}^{*}$  in (\ref{verticiQ*2}) is transformed into the unit square $Q$ in 10 steps.}
	\label{fig16}
\vskip1truecm
\centerline{\includegraphics[height=4cm]{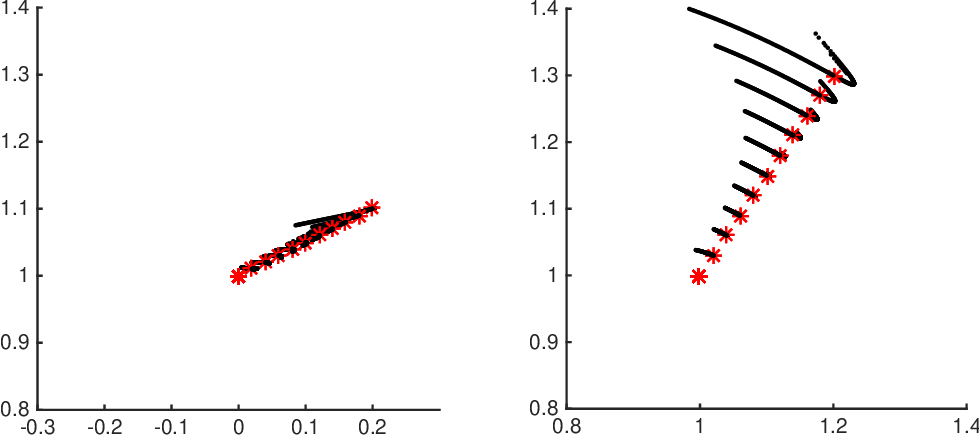} }
 \caption{Projections of the local curves joining the 
	corresponding isospectral domains for each of the quadrilaterals of Fig.~\ref{fig16}.
Unfortunately, the picture is not clear on the left-hand side. However, this test has
been reported in order to emphasize the interesting behavior developing on the right-hand side.}
 \label{fig17}
  \end{figure}

Finally, we run the same tests by using a different set of vertices $\hat V_k^*$, $k=1,2,3,4$, of the  initial quadrilateral  $ \hat{Q}^{*}$, namely:
 \begin{equation}\label{verticiQ*2}
 \hat{V}_{1}^{*}=(0,0), \ \ \hat{V}_{2}^{*}=(1,0), \ \ \hat{V}_{3}^{*}=(0.2 ,1.1), \ \ \hat{V}_{4}^{*}=(1.2 ,1.3).
  \end{equation}
In  Fig.~\ref{fig16} and Fig.~\ref{fig17} the reader can see the results of the tests obtained by 
taking  $S=10$ in (\ref{legge}) and  $T_0=0.06$ in  (\ref{legget}).  
Here the curve is twisting in a more complicated fashion.  Scattered
dots may appear in the plots. They are a consequence of the break down of the algorithm when
approaching the unit square $Q$. In order to get rid of them it is necessary to further
reduce the interval $[-T_j, T_j]$, as $j$ gets close to $S$.


%
\section{Conclusions}\label{sec8}

Still remaining in finite dimension, the extension of our analysis
to general domains depending on more degrees of freedom can be a severe
numerical task. The main difficulty is how to choose the family of admissible
domains. They may have, for instance, a polygonal boundary and the
approximation
of the continuous operator can be performed by finite elements. If we are
far from special symmetric configurations (as the case of the unit  square
considered in Sect.~\ref{sec7}), the local application of the Implicit Function theorem
should guarantee isospectral deformations in a small neighborhood.
The one-dimensional curve now belongs to a space of larger dimension and
its grafical representation can be rather troublesome. For visualization,
the best way
is probably to show animations concerning the movement of the isospectral
domains. Note that the costs of implementation can drastically become
high when
increasing the degrees of freedom.

An alternative to the work  presented  here is to try to preserve a certain
number of eigenvalues of the exact operator. For example, concerning the
family of quadrilaterals (as those examined so far), instead of introducing
a discretization based on matrices of dimension  $4\times 4$, one can evaluate the first
4 eigenvalues of the exact Laplacian and move the domains with the aim of
preserving their values. We expect the so obtained domains to be slightly
different
from the ones we found in this paper. Unfortunately, the exact eigenvalues of
$-\Delta$ on a general
quadrilateral domain are not explicitly available, so that the computation
should be accompanied by an appropriate discretization on a very fine
grid. Again, the complexity of the algorithm may become unaffordable
as more degrees of freedom are introduced.

The problem of finding isospectral families, connected with continuity,
in order to preserve
\underline{all} the eigenvalues of the \underline{exact} Laplace operator is certainly harder than the
experiments we tried in this paper, and represents a stimulating theoretical
challenge. We hope however that our little contribution may be the starting
point for future ideas.


\end{document}